\pdfoutput=1
\RequirePackage{ifpdf}
\ifpdf % We are running pdfTeX in pdf mode
\documentclass[pdftex]{sigma}
\else
\documentclass{sigma}
\fi

\begin{document}

%\allowdisplaybreaks

\renewcommand{\PaperNumber}{011}

\FirstPageHeading

\ShortArticleName{Deformation Quantization by Moyal Star-Product and Stratonovich Chaos}

\ArticleName{Deformation Quantization by Moyal Star-Product\\ and Stratonovich Chaos}

\Author{R\'emi L\'EANDRE~$^\dag$ and Maurice OBAME NGUEMA~$^\ddag$}

\AuthorNameForHeading{R.~L\'eandre and M.~Obame Nguema}

\Address{$^\dag$~Laboratoire de Math\'ematiques, Universit\'e de Franche-Comt\'e, 25030, Besan\c con, France}
\EmailD{\href{mailto:remi.leandre@univ-fcomte.fr}{remi.leandre@univ-fcomte.fr}}

\Address{$^\ddag$ Institut de Math\'ematiques de Bourgogne, Universit\'e de Bourgogne, 21000, Dijon, France}
\EmailD{\href{mailto:maurice-saint-clair.obame-nguema@u-bourgogne.fr}{maurice-saint-clair.obame-nguema@u-bourgogne.fr}}

\ArticleDates{Received November 16, 2011, in f\/inal form March 06, 2012; Published online March 15, 2012}

\Abstract{We make a deformation quantization by Moyal star-product on a space of functions endowed with the normalized Wick product and where Stratonovich chaos are well def\/ined.}


\Classification{60H40; 81S20}

\section{Introduction}

This work follows those of G.~Dito~\cite{GDi}  (motivated by quantum f\/ield theory
\cite{GDi2,GDi1,D-F,EW}), of R.~L\'eandre~\cite{R-S3,R-S,R-S1} and  G.~Dito, R.~L\'eandre~\cite{D-L}.
\cite{GDi}~deals with the deformation quantization of Moyal on a  Hilbert space:  the condition of equivalence of the  Moyal deformations is that the chosen perturbation operator is a Hilbert--Schmidt operator but in this case, Moyal and normal products are not equivalents.
 \cite{R-S} choose   Hida weighted Fock spaces  which are very small spaces. This gives spaces of continuous functionals on the path space.
 The inequivalences of~\cite{GDi} become equivalences. In~\cite{D-L}, the Malliavin
test algebra is used for  the Moyal quantization. A~very important remark  in~\cite{GDi} and~\cite{D-L} is that the matrix
 of the associated Poisson structure is still bounded.

Our motivation comes from the fact we deal with deformation
quantization on an algebra constituted of Stratonovich chaos. The Connes spaces where our work
was possible present some dif\/ferences with the Hida spaces of~\cite{R-S}. Indeed, Connes spaces are involved with tensor products of  Banach spaces and Hida spaces  are involved with tensor products of Hilbert spaces~\cite{HD1,R-L,NO}. In inf\/inite dimension analysis, there are two basic objects~\cite{PAM}:
\begin{itemize}\itemsep=0pt
\item an algebraic model;
\item a mapping space and a map $\Psi$ def\/ined from algebraic model into the space of functionals on this mapping space.
\end{itemize}\itemsep=0pt
Let us recall what is the main dif\/ference between the point of view of \cite{D-L} and the point of view of \cite{R-S,R-S1} and this paper:
\begin{itemize}\itemsep=0pt
\item  The tools of Malliavin calculus on the Wiener space are used in~\cite{D-L}. The Malliavin test algebra is constituted of functionals {\it almost surely defined}, because there is no Sobolev imbedding in inf\/inite dimension~\cite{MP1,MP2,DN,ASU}.
\item  The tools of white noise analysis are used  in~\cite{R-S,R-S1} and in this present work. The dif\/ferential operations and the topological structures are seen at the level of the algebraic model  and after they are {\it transported} through the map~$\Psi$  on a set of functionals which are {\it continuous} on the Wiener space. \cite{R-S}~and~\cite{R-S1} work on the level of the algebraic model. They use a dif\/ferent normalization of annihilation operators which f\/it with the map Stratonovitch chaos and not with the map Wiener chaos. In the present paper, we transport the algebraic operations of \cite{R-S,R-S1} on a functional space by using the map $\Psi$ Stratonovitch chaos which was not def\/ined in~\cite{R-S} and~\cite{R-S1}.
\end{itemize}
Note that E.~Getzler \cite{EG}  was the f\/irst to consider another map  than the map Wiener chaos. He used as algebraic space a Connes space and as map~$\Psi$ the map of Chen iterated integrals (see~\cite{R-S2} for developments).

After this quick  presentation, we shall def\/ine in the second part our Connes functional space
  on which  the multiple Stratonovich integrals or Stratonovich chaos are well def\/ined.
 We shall f\/inish this second part by a study of the annihilation operator on our  functions space.
 In the third part, the deformation quatization of the Poisson bracket by a Moyal product is def\/ined on the Connes space. The last part is about equivalences of  deformations on the Connes space.

\section{Gaussian space}

In this section we show the existence of a  Gaussian measure on a space of continuous loops which will be our reference measure throughout this work on the algebraic space. We also def\/ine Stratonovich chaos and dif\/ferentiation operators we use in the next section.

\subsection{Gaussian measure}

Let us consider the Hilbert space  $\mathcal{H}:=\mathcal{H}(S^1, \mathbb {R}^d)$
  such that $\gamma \in \mathcal{H}$  verif\/ies: $||\gamma ||^2 = \int_{0}^{1} |\gamma (s)|^{2} ds + \int_{0}^{1} |\dot\gamma (s)|^{2} ds $.
 We consider $B = \{B(t) = (B_i(t)),\, t \in S^{1} \}$ the Wiener process  associated to this Hilbert space. We  note $(\cdot, t) \mapsto G(\cdot, t)$
 the symmetric Green kernel. Let $h$ be a  continuous function  from $S^1$  with values in $\mathbb{R}$ such that
\begin{gather*}
h(1)=\int_{0}^{1}\langle h(s),G(s,1)\rangle ds+\int_{0}^{1}\langle \dot h(s),\frac{d}{ds}G(s,1)\rangle ds.
\end{gather*}
$G$ is solution of a second-class linear dif\/ferential equation and the  Green kernel  of that equation verif\/ies $\frac{d^2}{ds^2}G(s,1)-G(s,1)=0$ but also  $\frac{d}{ds}G(1,1)-\frac{d}{ds}G(0,1)=1$.
We obtain that $G(s,1)=\alpha e^{-s}+\beta e^{s}$ where $\alpha$ and $\beta$
 are constants of integration with respect to $s$.
 With the conditions in the limits, we f\/ind that $\alpha=\frac{-1}{2(1-e^{-1})}$ and $\beta=\frac{1}{2(1-e)}$. Moreover, we know that
 $E[B_i (s) B_j (t)] = \delta_{i,j}G(s,t)$ where $\delta_{i,j}$ denotes the classical Kronecker symbol. There is $\mu>0$ such that
$|G(t,t)+ G(s,s)-2G(s,t)|\leq \mu|t-s|$. By standard result on Gaussian measures, for all $ p>1$, there exist $\mu_p >0$ such that
\begin{gather*}
E\big[|B(t)-B(s)|^{2p}\big]\leq \mu_{p} |t-s|^p.
\end{gather*}
By the Kolmogorov's criterion of continuity, we deduct that $B$ is Hoelderian.

\subsection{Connes space}

Let us consider  the Hilbert space $\mathcal{H}$ above and a map $\gamma$ def\/ined from  the circle into $\mathbb{R}^{d}$ such that
\begin{gather*}
\int_{0}^{1}|\gamma(s)|^{2}ds + \int_{0}^{1}|\dot\gamma(s)|^{2}ds = ||\gamma||^{2}.
\end{gather*}

We consider a symmetric tensor product $F^n$
\begin{gather}\label{FN}
F^n = \frac{1}{n!}\sum_{\sigma\in \mathfrak{G}_{n}}F_{\sigma(1)}\otimes \cdots \otimes F_{\sigma(n)},
\end{gather}
where $\sigma$ is a permutation of the symmetric group of degree $n$. $F_i$ are elements of $\mathcal{H}$ and $\otimes$ denotes the standard
tensor product on this Hilbert space.

$F^n$ can be seen as a function from $(S^1)^n$ into $(\mathbb{R}^d)^{\otimes n}$. We consider the  symmetric Fock space constituted of $F= \sum F^n$
\begin{gather}\label{GN}
\sum (||F^{n}||^{\hat{{\otimes}}{n}})^2 < \infty,
\end{gather}
where we consider the Hilbert norm on the symmetric tensor product $\mathcal{H}^{\hat{\otimes}n}$.
We can def\/ine the Wiener chaos
\begin{gather*}
\langle F^{n},B(\cdot)\otimes \cdots \otimes B(\cdot) \rangle_W.
\end{gather*}
It is well def\/ined for the symmetric function $F^{n}\in \mathcal{H}^{\widehat{\otimes} n}$ and the Gaussian process $B$. The map Wiener chaos realizes as classical
an isomorphism between the symmetric Fock space and the $L^2$ of the Wiener space. The goal of this work is to replace Wiener chaos by Stratonovitch chaos.

 Let $\{e_{i}\}_{1\leq i\leq d}$ be a canonical basis of $\mathbb{R}^{d}$. We get by Fourier expansion  an orthonormal basis
 of the  Hilbert space  for some $\lambda >0 $
\begin{gather}\label{eq1}
\gamma_{i,k}(s)=\frac{\cos(2k \pi s)}{\sqrt{\lambda k^{2}+1}}e_{i}
\end{gather}
if $k\geq 0$ and if $k<0$
\begin{gather}\label{eq2}
\gamma_{i,k}(s)=\frac{\sin(2k \pi s)}{\sqrt{\lambda k^{2}+1}}e_{i}.
\end{gather}

Consider that $(\gamma_{i})_{i\geq 1}$ is an orthonormal basis of $\mathcal{H}$, then an orthonormal basis of $\mathcal{H}^{\widehat{\otimes} n}$
 takes the form
\begin{gather}\label{eq3}
\gamma_{N}(s)=\frac{1}{N!}\sum_{\sigma \in \frak{S}_{N}}\gamma_{\sigma(j_1)}(s_{\sigma(1)})\otimes\cdots\otimes\gamma_{\sigma(j_{N})}(s_{\sigma(N)}).
\end{gather}
$\otimes$  denotes the tensor product on~$\mathbb{R}^d$.
$\sigma$~is a permutation of the symmetric group $\frak{S}_{N}$ of degree~$N$.  For all  $l = 1,\ldots,N$
\begin{gather*}
j_{l} \in \mathfrak{J} = \bigcup_{p=1}^{m}\mathfrak{J}_{n_{p}}
\end{gather*}
and $\mathfrak{J}_{n_{p}} = \{(i_{p},k_{p}),\dots,(i_{p},k_{p})\}$ such that $|\mathfrak{J}_{n_{p}}| = n_{p}$ and $\mathfrak{J}_{n_{i}} \cap\mathfrak{J}_{n_{j}} = \varnothing$ if $i\neq j$ and we note by $|\mathfrak{J}| = N =\sum\limits_{p=1}^{m}n_{p}$ for $n_{p}\leq n$.

\begin{remark}\label{remark1}
We shall use the orthonormal basis of the symmetric space $\mathcal{H}^{\widehat{\otimes} n}$ given by \eqref{eq3} to avoid  redundancies in the calculations throughout this work.
\end{remark}

Let us consider  the space ${\rm CO}_{k,C}$ of function $F$ given by
\begin{gather*}
\sum F^{n} = F,
\end{gather*}
where  $F^{n} \in \mathcal{H}^{\widehat{\otimes} n}$ is  a $C^{\infty}$  function of $n$ parameters such that for all  $k, C>0$
\begin{gather*}
||F||_{k,C} = \sum C^{n}\sum_{\substack{J_{\alpha}\cap J_{\beta}=\varnothing \\ J_{1}\cup\cdots\cup J_{l} = \{1,\ldots,n\}}}\big\|D^{(n_{1})}_{S_{J_{1}}}\cdots D^{(n_{l})}_{S_{J_{l}}}F^{n}\big\|_{\infty}
 := \sum C^{n}||F^{n}||_{k} < \infty,
\end{gather*}
where  $J_{\alpha}=\{l_{1},\dots,l_{\alpha}\}$ with $J_{1}\cup\cdots\cup J_{l}=\{1,\ldots,n\}$ and $S_{J_{\alpha}}=\{s_{l_{1}},\ldots,s_{l_{\alpha}}\}$. For $\alpha\in \{1,\ldots,l\}$ such that  $n_{\alpha}\leq k$, we  write
\begin{gather*}
D^{(n_{\alpha})}_{S_{J_{\alpha}}}=\frac{\partial^{\sum n_{\alpha}}}{\partial s^{n_{\alpha}}_{l_{1}}\cdots \partial s^{n_{\alpha}}_{l_{\alpha}}}.
\end{gather*}

\begin{definition}\label{definition1}
We call Connes space the set   ${\rm CO}_{\infty-}$  given by
\begin{gather*}
{\rm CO}_{\infty-}=\bigcap_{k,C}{\rm CO}_{k,C}
\end{gather*}
for all $k,C > 0$.
\end{definition}

\begin{remark}\label{remark2}
A sequence $F_k$ of ${\rm CO}_{\infty-}$ converges to $F$ in ${\rm CO}_{\infty-}$ if it converges in all ${\rm CO}_{k,C}$.
\end{remark}
\begin{remark}\label{remark3}
Let us compare with the Hida Fock space. We consider the positive selfadjoint Laplacian $\Delta$  on $\mathcal{H}$  and we consider the $k^{\rm th}$ order Sobolev space~$\mathcal{H}_k$ associated with $\Delta+1$. It is  a Hilbert space. We consider the symmetric tensor product~$\mathcal{H}_k^{\widehat{\otimes}
n}$ endowed with its Hilbert norm~$\Vert\cdot \Vert_k^{\hat{n}}$. Let us consider a sequence~$F^n$ of $\mathcal{H}_k^{\widehat{\otimes}
n}$. $F = \sum F^n$ belongs to $W.N_{k,C}$ if $\sum C^n (\Vert F \Vert_k^{\hat{\otimes}n})^2$ is f\/inite. The Hida Fock space $W.N_{\infty-}$ is the intersection of all $W.N_{k,C}$. By Cauchy--Schwarz inequality and the Sobolev imbedding theorem we see that the Connes space is densely continuously imbedded in the Hida Fock space if we consider the standard Hilbert norm $\Vert\cdot\Vert_k^{\otimes n}$ on the symmetric tensor product of~$\mathcal{H}$ in the def\/inition of the Hida Fock space.
\end{remark}

\begin{definition}\label{definition2}
Let $F=\sum F^{n}$ and $G=\sum G^{m}$ be two functions of ${\rm CO}_{\infty-}$. We def\/ine Wick product of $F$ and $G$ by
\begin{gather}\label{tensor0}
{:F.G:} = \sum_{m,n}F^{n}\widehat{\otimes}G^{m}	.
\end{gather}
  $F^{n}\widehat{\otimes}G^{m}$ is the symmetric tensor product of the functions $F^{n}$ and $G^m$  given by
\begin{gather}\label{tensor}
F^{n}\widehat{\otimes}G^{m} = \frac{1}{(m+n)!}\sum_{\sigma\in \frak{S}_{n+m}} F^{n}\otimes_{\sigma} G^{m}	
\end{gather}
with $F^{n}$ and $G^m$ verifying both \eqref{FN} and \eqref{GN} and where  $\sigma$ is a permutation of the symmetric group  $\frak{S}_{n+m}$ of the space $\mathcal{H}^{\otimes(m+n)}$.
\end{definition}

\begin{theorem}\label{theorem1}
The Connes space  ${\rm CO}_{\infty-}$ is a topological commutative algebra for the  Wick pro\-duct.
\end{theorem}

\begin{proof}
For the  algebraic properties of the Wick product on ${\rm CO}_{\infty-}$, see Theorem~\ref{theorem2} and Proposition~\ref{proposition2}. Let us consider
\begin{gather*}
F=\sum F^{n} , \qquad G=\sum G^{m}.
\end{gather*}

By \eqref{tensor0} and  \eqref{tensor} we have
\begin{gather*}
{:F.G:} =\sum_{m,n}\frac{1}{(m+n)!}\sum_{\sigma\in \{1,\ldots,n+m\}} F^{n}\otimes_{\sigma} G^{m}.	
\end{gather*}

Clearly, for all $k > 0$, there exist $k_{0}>0$ such that
\begin{gather*}
||F^{n}\otimes_{\sigma} G^{m}||_{k}\leq ||F^{n}||_{k_{0}}||G^{m}||_{k_{0}}.
\end{gather*}
For $C > 0$, there exist $C_{0}>0$ such that
\begin{gather*}
||:F.G:||_{k,C}\leq \sum_{n,m}C^{m}_{0}C^{n}_{0}||F^{n}||_{k_{0}}||G^{m}||_{k_{0}}.
\end{gather*}
We deduce
\begin{gather*}
||:F.G:||_{k,C}\leq K||F||_{k_{0},C_{0}}||G||_{k_{0},C_{0}}.
\end{gather*}
 The theorem is proved.
\end{proof}

\subsection{Multiple Stratonovich integrals}

The theory of Stratonovich chaos was initiated in \cite{Hu-Me, So-U} but no convenient functional space was def\/ined.
Let us consider $\gamma_{N}=\gamma^{\otimes n_{1}}_{i_{1},k_{i_{1}}}\otimes\cdots\otimes\gamma^{\otimes n_{m}}_{i_{m},k_{i_{m}}}$ such that $\gamma^{\otimes n_{j}}_{i_{j},k_{i_{j}}}$ is the $n_{j}$-th tensor product of $\gamma_{i_{j},k_{i_{j}}} \in H$. The associated Stratonovitch chaos  takes the form
\begin{gather}\label{add-16}
I^{S}_{m}(\gamma_{N})= \prod_{1\leq j\leq m} \langle \gamma^{\otimes n_{j}}_{i_{j},k_{i_{j}}},B(\cdot)\otimes\cdots\otimes B(\cdot)\rangle_S,
\\
\label{add-17}
I^{S}_{m}(\gamma_{N})= \prod_{1\leq j\leq m} \langle \gamma_{i_{j},k_{i_{j}}}, B(\cdot)\rangle^{n_{j}}.
\end{gather}

We consider multiple Stratonovitch integrals. In Stratonovitch calculus, the It\^o formula reduces to the classical one.
 So \eqref{add-16} gives \eqref{add-17}, because in such a case the Stratonovitch--It\^o formula is
nothing else than  the ordinary Fubini theorem.  We have
\begin{gather*}
I^{S}_{m}(\gamma_{N})=\prod_{1\leq j \leq m}\bigg(\int_{S^1}\langle \gamma_{i_{j},k_{i_{j}}}(s_{j}),B(s_{j})\rangle ds_{j}+
\int_{S^1}\langle \dot\gamma_{i_{j},k_{i_{j}}}(s_{j}),\circ dB(s_{j})\rangle\bigg)^{n_{j}},
\end{gather*}
where we consider a Stratonovitch dif\/ferential $\circ dB(s_j)$.
With the integration by parts formula
\begin{gather*}
\int_{S^1}\langle \dot\gamma(s_{j}),dB(s_{j})\rangle= - \int_{S^1}\langle \frac{\partial^2}{\partial s^{2}_{j}}\gamma (s_{j}),B(s_{j})\rangle ds_{j},
\end{gather*}
we get
\begin{gather*}
I^{S}_{m}(\gamma_{N})=\prod_{1\leq j \leq m}\bigg(\int_{S^1}\langle \gamma_{i_{j},k_{i_{j}}}(s_{j}),B(s_{j})\rangle ds_{j}-
\int_{S^1}\langle \frac{\partial^2}{\partial s^{2}_{j}}\gamma_{i_{j},k_{i_{j}}}(s_{j}),B(s_{j})\rangle ds_{j}\bigg)^{n_{j}}.
\end{gather*}
Then
\begin{gather*}
I^{S}_{m}(\gamma_{N})=\int_{(S^1)^{N}}\prod_{l_{1}=1}^{n_{1}}\cdots \prod_{l_{n_{m}}=n_{m-1}+1}^{n_{m}}\big(\langle \gamma_{i_{1},k_{i_{1}}}(s_{l_{1}}),B(s_{l_{1}})\rangle \\
\phantom{I^{S}_{m}(\gamma_{N})=}{}
- \langle \frac{\partial^2}{\partial s^{2}_{l_{1}}}\gamma _{i_{1},k_{i_{1}}}(s_{1_{1}}),B(s_{l_{1}})\rangle\big)\cdots \big(\langle \gamma_{i_{n_{m}},k_{i_{n_{m}}}}(s_{l_{n_{m}}}),B(s_{l_{n_{m}}})\rangle \\
\phantom{I^{S}_{m}(\gamma_{N})=}{}
- \langle \frac{\partial^2}{\partial s^{2}_{l_{n_{m}}}}\gamma _{i_{n_{m}},k_{i_{n_{m}}}}(s_{l_{n_{m}}}),B(s_{l_{n_{m}}})\rangle \big) ds_{l_{1}}\cdots ds_{N},
\end{gather*}
and f\/inally
\begin{gather*}
I^{S}_{m}(\gamma_{N})=\sum_{J\subset\{1,\dots,N\}}\int_{(S^1)^{N}}(-1)^{N-|J|}\prod_{\substack{j\in J\\ l\notin J }}\langle \gamma_{i_{j},k_{i_{j}}}(s_{j}),B(s_{j})\rangle \langle \frac{\partial^2}{\partial s^{2}_{l}}\gamma_{i_{l},k_{i_{l}}}(s_{l}),B(s_{l})\rangle ds_{j}ds_{l}.
\end{gather*}

These considerations, which take care of the dif\/ference of behaviour between  Stratonovitch chaos and Wiener chaos into the passage from \eqref{add-16} to \eqref{add-17} allow us to give the next
def\/inition and
to  generalize by linearity  $I^{S}_{m}(\gamma_{N})$ in the functions  $F:=\sum F^{n}$ which are not products. We get

\begin{definition}\label{definition3}
The multiple Stratonovich  integrals or Stratonovich chaos takes the  form
\begin{gather*}
 I^{S}_{m}(F)=
\sum_{n}\sum_{J\subset \{1,\ldots,n\}}\int_{(S^{1})^{n}}(-1)^{n-|J|}\langle D^{(2)}_{S_{J}}F^{n}(s_{1},\ldots,s_{n}),
B(s_{1})\otimes\cdots \otimes B(s_{n})\rangle ds_{1}\cdots ds_{n},
\end{gather*}
and we put
\begin{gather*}
I^{S}_{m}(F):= \sum_{n}\langle F^{n}, B\otimes\cdots\otimes B\rangle_{S}.	
\end{gather*}
\end{definition}

\begin{remark}\label{remark4}
For arbitrary small  $\mu >0$, we can choose $n_{\mu}>0$($n_{\mu}$ depending only on $\mu$)
 such that for all $n > n_{\mu}$, we have  $\sup_{s}|D^{(2)}_{S_{J}}F^{n}(s)|\leq \mu^n$ and for all $\mathbf{M}>0$ such that  $||B||_{\infty}:=\sup(\{|B(s)|, s\in S^{1}\})\leq \mathbf{M} $, we  get
\begin{gather*}%\label{eq5}
\sup \big\{|\langle D^{(2)}_{S_{J}}F^{n}(s_{1},\ldots,s_{n}), B(s_{1})\otimes\cdots \otimes B(s_{n})\rangle\,|\,
 s_{1},\ldots,s_{n}\in S^1\big\} \leq \mathbf{M}^n\mu^{n}.
\end{gather*}
Then, there is  $C_{\mu}>0$  such that
\begin{gather*}%\label{eq6}
\sum _{m\geq 1}|| I^{S}_{m}(F)||_{\infty} \leq \sum_{n>n_{\mu}}2^{n}\mu^{n}\mathbf{M}^{n} + C_{\mu}<\infty,
\end{gather*}
where we def\/ine $|| I^{S}_{m}(F)||_{\infty}:= \sup(\{|I^{S}_{m}(F)|,||B||_{\infty}\leq \mathbf{M}\})$. Then $I^{S}_{m}(F)$ converges normally. Clearly  the map $ B \mapsto I^{S}_{m}(F)$  is continuous.
\end{remark}

\begin{theorem}\label{theorem2}
The multiple Stratonovich integrals  of the Wick product of two functions of ${\rm CO}_{\infty-}$ is equal  to the product of their multiple Stratonovich integrals  for all $F,G \in {\rm CO}_{\infty-}$
\begin{gather*}%\label{eq7}
I^{S}_{m}(:F.G:)=I^{S}_{m}(F).I^{S}_{m}(G).
\end{gather*}
\end{theorem}

\begin{proof}
We consider $F=\sum F^{n_{1}}$ and $G=\sum G^{n_{2}}$ two functions of ${\rm CO}_{\infty-}$. Since
\begin{gather*}
{: F.G :} = \frac{1}{(n_{1}+n_{2})!}\sum_{\sigma\in \frak{S}_{n_{1}+n_{2}}} F^{n_{1}}\otimes_{\sigma} G^{n_{2}},
\end{gather*}
we have
\begin{gather*}
I^{S}_{m}(:F.G:)= \sum_{n_{1},n_{2}}I^{S}_{m}\left( \frac{1}{(n_{1}+n_{2})!}\sum_{\sigma\in \frak{S}_{n_{1}+n_{2}}} F^{n_{1}}\otimes_{\sigma} G^{n_{2}}\right)\\
\phantom{I^{S}_{m}(:F.G:)}{} = \sum_{n_{1},n_{2}}\sum_{\sigma\in \frak{S}_{n_{1}+n_{2}}}\frac{1}{(n_{1}+n_{2})!}\langle F^{n_{1}}
\otimes_{\sigma} G^{n_{2}},B\otimes\cdots\otimes B\rangle_S.
\end{gather*}

But, by permutating indexes, we have clearly
\begin{gather*}
\langle F^{n_{1}}\otimes_{\sigma} G^{n_{2}},B\otimes\cdots\otimes B\rangle = \langle F^{n_{1}}\otimes G^{n_{2}},B\otimes\cdots\otimes B\rangle_{S}\\
\phantom{\langle F^{n_{1}}\otimes_{\sigma} G^{n_{2}},B\otimes\cdots\otimes B\rangle =}{}
\times \langle F^{n_{1}},B\otimes\cdots\otimes B\rangle_{S}\langle G^{n_{2}},B\otimes\cdots\otimes B\rangle_{S}.
\end{gather*}
By using  Fubini's theorem, we get
\begin{gather*}
I^{S}_{m}(:F.G:)=\sum_{n_{1},n_{2}}\langle F^{n_{1}},B\otimes\cdots\otimes B\rangle_{S}\langle G^{n_{2}},B\otimes\cdots\otimes B\rangle_{S}\\
\phantom{I^{S}_{m}(:F.G:)=}{}
=\left(\sum_{n_{1}}\langle F^{n_{1}},B\otimes\cdots\otimes B\rangle_{S}\right)
\left(\sum_{n_{2}}\langle G^{n_{2}},B\otimes\cdots\otimes B\rangle_{S}\right) = I^{S}_{m}(F).I^{S}_{m}(G).
\end{gather*}
 The theorem is proved.
\end{proof}

\subsection{Dif\/ferentiation operators}

Dif\/ferentiation operators are  annihilation and creation operators. In the case of Hida test algebra~\eqref{tensor},
these operators are adjoint operators and then their study is simplif\/ied. In our case, using Banach spaces to  def\/ine  the Connes space makes that it is dif\/f\/icult to give a def\/inition of an adjoint operator. Then, we just give  a description of annihilation operator.

\begin{definition}\label{definition4}
We def\/ine annihilation operator on ${\rm CO}_{\infty-}$, for all $ h\in \mathcal{H}$ and $F=\sum F^n\in {\rm CO}_{\infty-}$ by
\begin{gather*}
a_{h}F :=\sum \!\int_{S^1}\!\Big[\langle F^{n}(s_{1},\ldots,\bar s_{i},\ldots,s_{n-1}),h(\bar s_{i})\rangle
+\langle \frac{d}{d\bar s_{i}}F^{n}(s_{1},\ldots,\bar s_{i},\ldots,s_{n-1}),\dot h(\bar s_{i})\rangle\Big] d\bar s_{i},%\label{eq11}
\end{gather*}
where $\bar s_{i}$ means that we make a concatenation at this term.
\end{definition}

 We have
\begin{proposition}\label{proposition1}
The G\^ateaux derivative of a multiple Stratonovich integrals of a function $F\in {\rm CO}_{\infty-}$ is equal to the multiple Stratonovich integral of the annihilation operator of that function. For all  $ F\in {\rm CO}_{\infty-}$
\begin{gather}\label{eq8}
D_{h} I^{S}_{m}(F)=I^{S}_{m}(a_{h} F).
\end{gather}
\end{proposition}

\begin{proof}
We have
\begin{gather*}
D_{h}I^{S}_{m}(F)=
\sum_{n}\sum_i\sum_{J\subset \{1,\ldots,n\}}\int_{(S^1)^{n}}(-1)^{n-|J|}
\sum\langle D^{(2)}_{S_{J}}F^{n}(s_{1},\ldots,s_{i-1},\bar{s_{i}},s_{i+1},\ldots,s_{n-1}),\\
\hphantom{D_{h}I^{S}_{m}(F)=
\sum_{n}\sum_i\sum_{J\subset \{1,\ldots,n\}}}{}
B(s_{1})\otimes\cdots\otimes h(\bar s_{i})\otimes\cdots\otimes B(s_{n})\rangle ds_{1}\cdots ds_{n}\\
\hphantom{D_{h}I^{S}_{m}(F)}{}
=\sum_{n}\sum_{J\subset \{1,\ldots,n-1\}}\int_{(S^1)^{n-1}}(-1)^{n-|J|-1}
\langle D^{(2)}_{S_{J}}a_{h}F^{n}(s_{1},\ldots,s_{n-1}),\\
\hphantom{D_{h}I^{S}_{m}(F)=
\sum_{n}\sum_i\sum_{J\subset \{1,\ldots,n\}}}{}
B(s_{1})\otimes\cdots\otimes B(s_{n-1})\rangle ds_{1}\cdots ds_{n-1}.
\end{gather*}
The result holds.
\end{proof}

\begin{remark}\label{remark5}
Let be $h \in \mathcal{H}$. There is $C_1> C$ such that for  for all $F=\sum F^n \in {\rm CO}_{\infty-}$
\begin{gather*}
||a_{h} F||_{k,C}\leq \sum C^{n}_{1}||h|| ||F^{n+1}||_{k}<\infty.
\end{gather*}
 The annihilation operator is continuous on ${\rm CO}_{\infty-}$.
\end{remark}

\begin{proposition}\label{proposition2}
The annihilation operator is a derivation for the Wick product on ${\rm CO}_{\infty-}$.  For all $F,G\in {\rm CO}_{\infty-}$
\begin{gather*}
a_{h}(: F.G :) =  : (a_{h} F) . G : + : F . (a_{h} G ):.
\end{gather*}
\end{proposition}

\begin{proof}
We have just to  show that the map  $F\mapsto I^{S}_{m}(F)$ is injective. Indeed, considering~\eqref{eq8} and the fact that  $a_{h}$ is a derivation  on the algebraic space , it is clear that if $F \mapsto I^{S}_{m}(F)$ is injective, the proposition is proved.

We suppose that for $F\in {\rm CO}_{\infty-} $, we have $I^{S}_{m}(F) = 0$ and for $z\in \mathbb{C}$, we put $\phi(z)=I^{S}_{m}(z.F)$. We get the following power series
\begin{gather*}
\phi(z)=\sum_{n}\!\! \sum_{J\subset \{1,\ldots,n\}}\! \int_{(S^{1})^{n}}\!\!(-1)^{n-|J|}\langle D^{(2)}_{S_{J}}F^{n}(s_{1},\ldots,s_{n}),
 z.B(s_{1})\otimes\cdots \otimes z.B(s_{n})\rangle ds_{1}\cdots ds_{n}\\
\phantom{\phi(z)}{}
=\sum_{n}z^{n}\!\!\sum_{J\subset \{1,\ldots,n\}}\!\int_{(S^{1})^{n}}\!\!(-1)^{n-|J|}\langle D^{(2)}_{S_{J}}F^{n}(s_{1},\ldots,s_{n}), B(s_{1})
\otimes\cdots \otimes B(s_{n})\rangle ds_{1}\cdots ds_{n}.
\end{gather*}

Since $I^{S}_{m}(F)= 0 $, we have $\phi(z)=0$. We deduce that  for all  $  n\geq 0$
\begin{gather*}
\sum_{J\subset \{1,\ldots,n\}}\int_{(S^{1})^{n}}(-1)^{n-|J|}\langle D^{(2)}_{S_{J}}F^{n}(s_{1},\ldots,s_{n}), B(s_{1})
\otimes\cdots \otimes B(s_{n})\rangle ds_{1}\cdots ds_{n}=0
\end{gather*}
for all $ h\in \mathcal{H}$. It can be written that
\begin{gather*}
\sum_{J\subset \{1,\ldots,n\}}\int_{(S^{1})^{n}}(-1)^{n-|J|}\langle D^{(2)}_{S_{J}}F^{n}(s_{1},\ldots,s_{n}),  h(s_{1})
\otimes\cdots \otimes h(s_{n})\rangle ds_{1}\cdots ds_{n}=0.
\end{gather*}

  We use Meyer's isomorphism  from $\mathcal{C}^{\infty}(S^{1},\mathbb{R}^{d})$ into $\mathcal{C}^{\infty}([0,1],\mathbb{R}^{d})$ given by
\begin{gather*}
\gamma := \sum_{i,k} a_{i,k}\gamma_{i,k}\mapsto \sum_{i}\sum_{k<0}a_{i,k}\alpha_{k}\gamma_{i,k} + \sum_{i} a_{i,o}\tilde{\gamma}_{i,0} + \sum_{i}\sum_{k\geq 0}a_{i,k+1}\tilde{\alpha}_{k}\gamma_{i,k}.
\end{gather*}
The coef\/f\/icients $a_{i,k}$ are for fast diminution, $\tilde{\gamma}_{i,0}(s)= s.e_{i}$ with $|\alpha_{k}|<\infty$ for $k>0$ and $|\tilde{\alpha}_{k}|<\infty$ for $k<0$. For all $\gamma$ we have $\int_{0}^{1}|\gamma(s)|^{2}ds<\infty$ with $\gamma(0)=0$. We deduce
\begin{gather*}
\int_{0<s_{1}<\cdots<s_{n}<1}F(s_{1},\ldots,s_{n})h(s_{1})\cdots h(s_{n})ds_{1}\cdots ds_{n} = 0
\end{gather*}
for all $h \in \mathcal{C}^{\infty}(S^{1},\mathbb{R}^{d})$,  where $F$ is smooth symmetric.
$F$ is a smooth symmetric function
\begin{gather*}
F(s_{1},\ldots,s_{n})=\sum_{i_{1},\ldots,i_{d}\geq 1}F^{i_{1}\cdots i_{d}}(s_{1},\ldots,s_{n})e_{i_{1}}\otimes\cdots\otimes e_{i_{d}}.
\end{gather*}
Let us suppose that $F\not= 0$.
Without restriction we can suppose that there exists an $\epsilon>0$, there  exist
$s^{(\varepsilon)}_{1}< \cdots<s^{(\varepsilon)}_{d}$ and $i_{1},\dots, i_{d}$ such that
\begin{gather*}
F^{i_{1}\cdots i_{d}}\big(s^{(\varepsilon)}_{1},\ldots,s^{(\varepsilon)}_{d}\big)>\varepsilon.
\end{gather*}

Since $F$ is smooth, there is very small $\eta > 0$ such that
\begin{gather*}
F^{i_{1}\cdots i_{d}}\big(s^{(\varepsilon)}_{1},\ldots,s^{(\varepsilon)}_{n}\big)>\frac{\varepsilon}{2}
\end{gather*}
on the product $[s^{(\varepsilon)}_{1}-\eta,s^{(\varepsilon)}_{1}]\times \cdots\times [s^{(\varepsilon)}_{d}-\eta,s^{(\varepsilon)}_{d}] = \prod I_{k}$. Then we take set $h =\sum 1_{I_{k}}e_{i_{k}}$ which give a contradiction. This shows  if $I^{S}_{m}(F)=0$, we have necessary $F = 0$. Thus, the map $F \mapsto I^{S}_{m}(F)$ is injective and the proposition is proved.
\end{proof}

\section{Poisson space}

The theory of deformation quantization was initiated in \cite{BFFS,BFFS2}. See \cite{DS,MY,AW} for reviews and~\cite{MG} for basical background.
This section gives some properties of the Poisson structure of the Connes space ${\rm CO}_{\infty-}$. We make also the quantization deformation of that Poisson structure in  Moyal star-product. We note by $\mathbb{K}= \mathbb{R}$ or $\mathbb{C}$.

\begin{definition}\label{definition5}
A Poisson structure on ${\rm CO}_{\infty-}$ is given by a $\mathbb{K}$-bilinear map $\{\cdot,\cdot \}$ from ${\rm CO}_{\infty-}\times {\rm CO}_{\infty-}$ into  ${\rm CO}_{\infty-}$ such that:
\begin{enumerate}\itemsep=0pt
\item $\{\cdot,\cdot \}$ is antisymmetric, satisf\/ies the Jacobi identity  and verif\/ies the Leibniz rule for the Wick product of ${\rm CO}_{\infty-}$.
\item For all $k,C$, there exists $K, k_1, C_1$ such that for all $F_1, F_2  \in {\rm CO}_{\infty-}$ we get
\begin{gather*}
||\{F_{1},F_{2}\}||_{k,C}\leq K||F_{1}||_{k_{1},C_{1}}||F_{2}||_{k_{1},C_{1}}.
\end{gather*}
\end{enumerate}
\end{definition}
We note by ${\rm CO}_{\infty-}[[\hbar]]$ the set of formal series with coef\/f\/icients in the Connes space ${\rm CO}_{\infty-}$.

\begin{definition}\label{definition6}
A star-product on ${\rm CO}_{\infty-}[[\hbar]]$ is a  $\mathbb{K}[[h]]$-bilinear map $\star_{\hbar}$ on ${\rm CO}_{\infty-}[[\hbar]]\times {\rm CO}_{\infty-}[[\hbar]]\!$ valued in ${\rm CO}_{\infty-}[[\hbar]]$. For all $ F_{1},F_{2}\in {\rm CO}_{\infty-}$ we have
\begin{gather*}
F_{1}\star_{\hbar}F_{2} = \sum_{r\geq 0}\frac{\hbar^{r}}{r!}P^{r}(F_{1},F_{2}).
\end{gather*}

For all $r\geq 0$, $P^{r}$ is a bilinear map on ${\rm CO}_{\infty-}$ satisfying:
\begin{enumerate}\itemsep=0pt
\item $P^{0}(F_{1},F_{2})=:F_{1}.F_{2}:.$
\item $P^{1}(F_{1},F_{2})-P^{1}(F_{2},F_{1})= 2\{F_{1},F_{2}\}.$
\item For all $r >0$, for all $  k , C >0$ there are $K,k_{1},C_{1}>0$ such that for all $F_1, F_2 \in {\rm CO}_{\infty-}$, we get
\begin{gather*}
||P^{r}(F_{1},F_{2})||_{k,C}\leq K||F_{1}||_{k_{1},C_{1}}||F_{2}||_{k_{1},C_{1}}.
\end{gather*}

\item For all $F_{1},F_{2},F_{3}\in {\rm CO}_{\infty-}$, we have: $F_{1}\star_{\hbar}(F_{2}\star_{\hbar}F_{3})=(F_{1}\star_{\hbar}F_{2})\star_{\hbar}F_{3}$.
\end{enumerate}
\end{definition}

We call $\star_{\hbar}$ a deformation of the Poisson bracket on ${\rm CO}_{\infty-}$.
\begin{definition}\label{definition7}
Two deformation quantizations $\star^{1}_{\hbar}$ and  $\star^{2}_{\hbar}$ of the same Poisson bracket are said equivalent if there exists a $\mathbb{K}[[\hbar]]$-linear map $T : {\rm CO}_{\infty-}[[\hbar]] \to {\rm CO}_{\infty-}[[\hbar]]$ expressed as a formal $T = I + \sum_{r\geq 1}\hbar^{r}T_{r}$  satisfying:
\begin{enumerate}\itemsep=0pt
\item  For all  $r\geq 1$,  $T_{r}: {\rm CO}_{\infty-}\to {\rm CO}_{\infty-} $  is a continuous operators and $T_{0}$ is the identity operator of ${\rm CO}_{\infty-}$.
\item  For all $ F,G\in {\rm CO}_{\infty-}$, we have $T(F\star^{1}_{\hbar}G)= T(F)\star^{2}_{\hbar}T(G).$
\end{enumerate}
\end{definition}

Let us consider  $\omega =(\omega_{ij})_{i,j\geq 1}$ a non-degenerate bilinear and antisymmetric form in  $\mathbb{R}^{d}\oplus\mathbb{R}^{d}$. We def\/ine the bilinear and antisymmetric form $\Omega$ on  $\mathcal{H}(S^1,\mathbb{R}^{d})$.  For all $\gamma_{1} , \gamma_{2} \in \mathcal{H} $
\begin{gather*}
\Omega(\gamma_{1},\gamma_{2})=\int_{S^1}\omega(\gamma_{1}(s),\gamma_{2}(s))ds.
\end{gather*}

We recall that  $\gamma_{i,k}$ are given by  \eqref{eq1} and \eqref{eq2}. We have
\begin{gather*}
\Omega(\gamma_{i,k_{i}},\gamma_{j,k_{j}})=\frac{\omega_{ij}}{Ck^2+1} \delta_{k_{i}k_{j}},
\end{gather*}
where $\delta_{k_{i}k_{j}}$ is the Kronecker delta. We note by  $\omega_{ij} = (\omega^{ij})^{-1}$. The Poisson matrix of $\Omega$ is given~by
\begin{gather*}
\{\gamma_{i,k_{i}},\gamma_{j,k_{j}}\}= \big(Ck^2+1\big). \omega^{ij} . \delta_{k_{i}k_{j}}.
\end{gather*}

Let us note by $a_{\gamma_{i,k}}$ the annihilation operator associated to $\gamma_{i,k}$. Then, for all $F,G\in {\rm CO}_{\infty-}$ we have
\begin{gather}\label{eq6*}
\{F , G\}=\sum_{m,n}\sum_{i,j,k} \big(Ck^2+1\big) \omega^{ij}\times :a_{\gamma_{i,k}}F^{n}.a_{\gamma_{j,k}}G^{m}:.
\end{gather}

\begin{proposition}\label{proposition3}
$\{\cdot,\cdot \}$ defines  a Poisson structure  on ${\rm CO}_{\infty-}$ in the sense of Definition~{\rm \ref{definition5}}.
\end{proposition}

\begin{proof}
Since the Wick product is continuous, for all $k,C>0$, there are $k_{2} , C_{2}>0 $ such that
\begin{gather*}
||\{F , G\}||_{k,C}\leq \sum_{m,n,k}\big(Ck^{2}+1\big)||a_{\gamma_{i,k}}F^{n}||_{k_{2},C_{2}}||a_{\gamma_{j,k}}G^{m}||_{k_{2},C_{2}}.
\end{gather*}
Using integration by parts, there are $k_{3} , C_{3}>0$ large enough such that
\begin{gather*}
||\{F,G\}||_{k,C}\leq \sum_{m,n,k}\big(Ck^{2}+1\big)^{-1}||F^{n}||_{k_{3},C_{3}}||G^{m}||_{k_{3},C_{3}}.
\end{gather*}
Then
\begin{gather*}
||\{F , G\}||_{k,C}\leq \left(\sum||F^{n}||_{k_{3},C_{3}}\right)\left(\sum||G^{n}||_{k_{3},C_{3}}\right)\left(\sum \frac{1}{(Ck^2+1)}\right).
\end{gather*}
Finally
\begin{gather*}
||\{F , G\}||_{k,C}\leq K||F||_{k_{3},C_{3}}||G||_{k_{3},C_{3}}< \infty.\tag*{\qed}
\end{gather*}
  \renewcommand{\qed}{}
\end{proof}

\begin{remark}\label{remark6}
Integration by parts  allows us to change the factor $(Ck^2+1)$ into $(Ck^{2}+1)^{-1}$ in the proof and give a bounded form of the Poisson bracket. Thus,  $\Omega$ acts continuously  on the space~${\rm CO}_{\infty-}$.
\end{remark}

Using the Wick product,  we can def\/ine the powers of the Poisson bracket as following: for all $ r\geq 0$ and $F,G \in {\rm CO}_{\infty-}$
\begin{gather*}
P^{r}(F,G)=\sum_{n,m}\!\sum_{\substack{i_{1},\ldots,i_{r}\geq 1\\j_{1},\ldots,j_{r}\geq 1\\k_{1},\ldots,k_{r}\geq 1}}\!\prod_{j=1}^{r}\!\big(Ck^{2}_{j}+1\big)\omega^{i_{1}j_{1}}\cdots \omega^{i_{r}j_{r}}
%\times
 :a_{\gamma_{i_{1},k_{1}}}\!\cdots a_{\gamma_{i_{r},k_{r}}}F^{n}.a_{\gamma_{j_{1},k_{1}}}\!\cdots a_{\gamma_{j_{r},k_{r}}}G^{m}.
\end{gather*}
Then we have
\begin{definition}\label{definition8}
The Moyal star-product on  ${\rm CO}_{\infty-}$ is given by
\begin{gather}\label{eq9}
F\star_{\hbar}G= :F.G: + \sum_{r\geq 1}\frac{\hbar^{r}}{r!}P^{r}(F,G).
\end{gather}
\end{definition}

The  Moyal star-product endowed with  the symplectic structure of $\Omega$ is well def\/ined on ${\rm CO}_{\infty-}$ in the
sense of Def\/inition~\ref{definition6}. Then, we have

\begin{theorem}\label{theorem3}
The formula \eqref{eq9} defines a deformation quantization of $\{\cdot,\cdot \}$ on ${\rm CO}_{\infty-}$ in the framework of Definition~{\rm \ref{definition6}}.
\end{theorem}
\begin{proof}
Since the Wick product is continuous, for all $k,C>0$ there are $k'_{1},C'_{1}>0$ such that
\begin{gather*}
||P^{r}(F,G)||_{k,C} \leq \sum_{m,n,k}\prod_{j}(Ck^{2}_{j}+1)||a_{\gamma_{i_{1},k_{1}}}\cdots a_{\gamma_{i_{r},k_{r}}}F^{n}||_{k'_{1},C'_{1}}
 ||a_{\gamma_{j_{1},k_{1}}}\cdots a_{\gamma_{j_{r},k_{r}}}G^{m}||_{k'_{1},C'_{1}}.
\end{gather*}
Using integration by parts, there are $k'_{2},C'_{2}>0$ large enough such that
\begin{gather*}
||P^{r}(F,G)||_{k,C}\leq \sum \prod_{j}^{r}(Ck^{2}_{j}+1)^{-1}||F^{n}||_{k'_{2},C'_{2}}||G^{m}||_{k'_{2},C'_{2}}.
\end{gather*}
Without loss of generality, we get
\begin{gather*}
||P^{r}(F,G)||_{k,C}\leq\left(\sum ||F^{n}||_{k'_{2},C'_{2}}\right)\left(\sum ||G^{m}||_{k'_{2},C'_{2}}
\right)\left(\sum \prod_{j}^{r} \frac{1}{(Ck^{2}_{j}+1)}\right).
\end{gather*}
Then
\begin{gather*}
||P^{r}(F,G)||_{k,C}\leq K'||F||_{k'_{2},C'_{2}}||G||_{k'_{2},C'_{2}}<\infty.
\end{gather*}

It remains to check the algebraic properties. It is enough to prove them if we consider f\/inite sum of~$\gamma_N$ because~$P^r$ apply
 the product of this space on itself and because by Stone--Weierstrass theorem the set of f\/inite sum of~$\gamma_N$ is dense in $C0_{\infty-}$. But in~\cite{R-S}, these algebraic properties were proved where a completion of Hida type of the set of f\/inite sum (by using an Hida Fock space) was chosen.

The result holds.
\end{proof}

\begin{remark}\label{remark7}
Since the map $F\mapsto I^{S}_{m}(F)$ is injective, we can use the dictionary between  the multiple Stratonovich integrals and the algebraic model  ${\rm CO}_{\infty-}$. Then,
the formula \eqref{eq9} becomes
\begin{gather*}
I^{S}_{m}(F)\star_{\hbar}I^{S}_{m}(G)= I^{S}_{m}(F).I^{S}_{m}(G) + \sum_{r\geq 1}\frac{\hbar^{r}}{r!}P^{r}(I^{S}_{m}(F),I^{S}_{m}(G)).
\end{gather*}
On $I^S_m{F}$ we choose the Banach norm of $F$.
\end{remark}

\section{Equivalence of deformation quantization}

Using the model of~\cite{GDi},we show that there are many similarities between the Connes space we use here  and the Hida test algebra of~\cite{R-S}.
Let us consider the Hilbert space $\mathcal{H} =\mathcal{H}(S^1,\mathbb{R}^{d})$ of functions $\gamma$ def\/ined from the circle into~$\mathbb{R}^{d}$ such that
\begin{gather*}
||\gamma||^{2} = \int_{0}^{1}|\gamma(s)|^{2}ds +\int_{0}^{1}\left|\frac{d}{ds}\gamma(s)\right|^{2}ds< \infty.
\end{gather*}

Let $(e_{i})_{1\leq i\leq d}$ be the canonical basis of $\mathbb{R}^{d}$. Considering the Fourier
 basis of $\mathcal{H}$ def\/ined by~\eqref{eq1} and~\eqref{eq2}, we can def\/ine on the Hilbert space $\mathcal{H}\oplus\mathcal{H}^{*}$ a Poisson
 structure by
 $\Omega(\Gamma_{1},\Gamma_{2})=\int_{0}^{1}\omega(\Gamma_{1}(s),\Gamma_{2}(s))ds +\int_{0}^{1}
\omega(\frac{d}{ds}\Gamma_{1}(s),\frac{d}{ds}\Gamma_{2}(s))ds$, where $\omega = (\omega_{ij})_{i,j\geq 1}$
 is a non-degenerated bilinear and antisymmetric form of $\mathbb{R}^{d}\oplus \mathbb{R}^{d}$ such that for all  $i\not= j$
\begin{gather*}
\omega_{ij} = 0, \qquad \omega_{ii^*} = 1, \qquad \omega_{i^*j^*} = 0, \qquad \omega_{i^*i} = -1,
\end{gather*}
and note $\Gamma_{j\in\{1,2\}}=\gamma_{j}\oplus\gamma^{*}_{j}\in \mathcal{H}\oplus\mathcal{H}^{*}$. We get that $\Omega$
 acts continuously on ${\rm CO}_{\infty-}\times {\rm CO}_{\infty-}$ and its Poisson  matrix is bounded. Then, the model of~\cite{GDi}  holds for the rest of the section.

\begin{definition}\label{definition9}
For all $\gamma\oplus \gamma^{*}\in\mathcal{H}\oplus\mathcal{H}^{*}$, we call Wick  exponentials,   the maps $\Phi_{\gamma,\gamma^{*}}$  def\/ined by
\begin{gather*}
h\oplus h^{*} \in \mathcal{H}\oplus\mathcal{H}^{*}\mapsto \Phi_{\gamma,\gamma^{*}}(h,h^{*}):=
 \exp(\langle h,\gamma\rangle+\langle h^{*},\gamma^{*}\rangle).
\end{gather*}
\end{definition}
We get a classical result for Hida calculus given by

\begin{proposition}\label{proposition4}
The Wick exponentials are  dense in ${\rm CO}_{\infty-}$.
\end{proposition}

\begin{proof}
We shall note by ${\rm CO}^{W}_{k,C}$ the adherence of Wick exponentials in ${\rm CO}_{k,C}$ and by ${\rm CO}^{n}_{k,C}$ the space of the
 product of $n$ homogeneous polynomials of ${\rm CO}_{k,C}$. We are going to use recurrence on the holomorphic function
\begin{gather*}
F(\lambda)=\exp\left[\lambda\left(\int_{S^{1}}\langle \Gamma(s),T(s)\rangle ds + \int_{S^{1}}\langle \frac{d}{ds}\Gamma(s),dT(s)\rangle\right)\right],
\end{gather*}
where $\Gamma=\gamma\oplus \gamma^{*}$ and $T= B\oplus B^{*}$. Then, $F(\lambda)$ can be written under the form
\begin{gather*}
F(\lambda)= \sum_{n\geq 0}\frac{\lambda^n}{n!}\left(\int_{S^{1}}\langle \Gamma(s),T(s)\rangle ds + \int_{S^{1}}\langle \frac{d}{ds}\Gamma(s),dT(s)\rangle\right)^n,
\end{gather*}
and obviously $F(\lambda)\in {\rm CO}^{W}_{k,C}$. With Cauchy's dif\/ferentiation formula
\begin{gather*}
F^{(n)}(\lambda)=\frac{n!}{2i\pi}\int_{S^1}\frac{F(z)}{(z-\lambda)^{n+1}}dz.
\end{gather*}
 It  is clear that  for all $n\geq 0$
\begin{gather*}
\left(\int_{S^{1}}\langle \Gamma(s),T(s)\rangle ds + \int_{S^{1}}\langle
 \frac{d}{ds}\Gamma(s),dT(s)\rangle\right)^{n}\in {\rm CO}^{W}_{k,C}.
\end{gather*}
Now, it remains just to prove that all products of $n$ homogeneous polynomials are in the adherence. We consider for  $|z|<1$ the holomorphic function
\begin{gather*}
F_{n+1}(z)= \int_{S^{1}}\langle\Gamma(s)+ z.\Gamma_{1}(s),T(s)\rangle ds +\int_{S^{1}}\left(\langle \frac{d}{ds}\Gamma(s) +
 z.\frac{d}{ds}\Gamma_{1}(s),dT(s)\rangle\right)^{n+1}.
\end{gather*}
$F_{n+1}$ is clearly in ${\rm CO}^{W}_{k,C}$ and by Cauchy's dif\/ferentiation formula
\begin{gather*}
F'_{n+1}(z)= \frac{1}{2i\pi}\int_{S^{1}}\frac{F(u)}{(u-z)^{2}}du.
\end{gather*}
Then $F'_{n+1}$ is also a function of ${\rm CO}^{W}_{k,C}$. By computation, we get
\begin{gather*}
F'_{n+1}(z)_=(n+1)\left(\int_{S^{1}}\langle\Gamma_{1}(s),T(s)\rangle ds+\int_{S^{1}}\frac{d}{ds}\langle\Gamma_{1}(s),dT(s)\rangle\right).F_{n}(z).
\end{gather*}
Then
\begin{gather*}
F'_{n+1}(0)=(n+1)\left(\int_{S^{1}}\langle\Gamma_{1}(s),T(s)\rangle ds+\int_{S^{1}}\frac{d}{ds}\langle\Gamma_{1}(s),dT(s)\rangle\right)\\
\phantom{F'_{n+1}(0)=}{} \times\left(\int_{S^{1}}\langle\Gamma(s),T(s)\rangle ds+\int_{S^{1}}\langle\frac{d}{ds}\Gamma(s),dT(s)\rangle\right)^n.
\end{gather*}

Thus, we proved the recurrence relation in the order $(n+1)$. We have
\begin{gather*}
\sum\prod_{i=1}^{n}\int_{S^{1}}\langle \Gamma_{i}(s),T(s)\rangle ds + \int_{S^{1}}\langle \frac{d}{ds}\Gamma_{i}(s),dT\rangle \in {\rm CO}^{W}_{k,C}.
\end{gather*}
By the theorem of Stone--Weierstrass, for all $ k , C >0$ we get that
\begin{gather*}
I^{S}_{m}(F) \in {\rm CO}^{W}_{k,C},
\end{gather*}
because $F$ is a limit of elements of ${\rm CO}^{n}_{k,C}$. We conclude that
\begin{gather*}
\bigcap_{k,C}{\rm CO}^{W}_{k,C}={\rm CO}^{W}_{\infty-} = {\rm CO}_{\infty-}.
\end{gather*}
 The proposition is proved.
 \end{proof}

According to \cite{R-S}, we choose the operator $A :\gamma_{i,k}\mapsto \alpha_{k}\gamma_{i,k}$ such that $|\alpha_{k}|\leq K|k|^{\mu}$  for some $\mu > 0$. We put
\begin{gather*}%\label{ea}
E_{A}\big(I^{S}_{m}(F),I^{S}_{m}(G)\big)=\sum_{i,k\geq 1}\big[\alpha_{k}D_{\gamma_{i,k}}I^{S}_{m}(F)D_{\gamma^{*}_{i,k}}I^{S}_{m}(G)
 +\alpha_{k}D_{\gamma_{i,k}}I^{S}_{m}(G)D_{\gamma^{*}_{i,k}}I^{S}_{m}(F)\big],
\end{gather*}
where $D_{\gamma_{i,k}}$(resp.\ $D_{\gamma^{*}_{i,k}}$) is the G\^ateaux derivative at $\gamma_{i,k}$ (resp.\ $\gamma^{*}_{i,k}$) in the direction $\mathcal{H}$ (resp.\ $\mathcal{H}^{*}\sim\mathcal{H}$).

\begin{theorem}\label{theorem4}
$E_{A}$ is continuous from ${\rm CO}_{\infty-}\times {\rm CO}_{\infty-}$ into ${\rm CO}_{\infty-}$.
\end{theorem}
\begin{proof}
For all $k, C >0$, integrating by parts we can f\/ind $k_{2},C_{2}>0$ large enough such that
\begin{gather*}
\big\|E_{A}\big(I^{S}_{m}(F),I^{S}_{m}(G)\big)\big\|_{k,C}\leq \sum \big(Ck^{2}+1\big)^{-1}\big\|I^{S}_{m}(F)\big\|_{k_{2},C_{2}} \big\|I^{S}_{m}(G)\big\|_{k_{2},C_{2}}<\infty.
\end{gather*}
This proves the theorem.
\end{proof}

We put
\begin{gather*}%\label{eqC}
C^{A}_{1}\big(I^{S}_{m}(F),I^{S}_{m}(G)\big)=\big\{I^{S}_{m}(F),I^{S}_{m}(G)\big\} + E_{A}\big(I^{S}_{m}(F),I^{S}_{m}(G)\big).
\end{gather*}
Using \eqref{eq6*}, we get
\begin{gather*}
C^{A}_{1}(I^{S}_{m}(F),I^{S}_{m}(G))=\sum_{i,k\geq 1}\big[(\alpha_{k}+1) D_{\gamma_{i,k}}I^{S}_{m}(F)D_{\gamma^{*}_{i,k}}I^{S}_{m}(G)
\\ \nonumber
\phantom{C^{A}_{1}(I^{S}_{m}(F),I^{S}_{m}(G))=}{}
 + (\alpha_{k}-1)D_{\gamma_{i,k}}I^{S}_{m}(G)D_{\gamma^{*}_{i,k}}I^{S}_{m}(F)\big].
\end{gather*}
Then,  we put  in the sense of dif\/ferential operators
\begin{gather*}
C^{A}_{r}\big(I^{S}_{m}(F),I^{S}_{m}(G)\big)=\big(C^{A}_{1}\big)^{r}\big(I^{S}_{m}(F),I^{S}_{m}(G)\big).
\end{gather*}
 We get
\begin{gather*}
C^{A}_{r}\big(I^{S}_{m}(F),I^{S}_{m}(G)\big)=\sum_{\substack{i_{1},\ldots,i_{r}\geq 1\\k_{1},\ldots,k_{r}\geq 1}}\prod_{l=1}^{r} (\alpha_{k_{l}}+1)D_{\gamma^{\sharp}_{i_{1},k_{1}}}\cdots D_{\gamma^{\sharp}_{i_{r},k_{r}}}I^{S}_{m}(F)\nonumber\\
\phantom{C^{A}_{r}\big(I^{S}_{m}(F),I^{S}_{m}(G)\big)=}{}
\times(\alpha_{k_{l}}-1)D_{\gamma^{\sharp}_{i_{1},k_{1}}}\cdots D_{\gamma^{\sharp}_{i_{r},k_{r}}}I^{S}_{m}(G).\label{eqC1}
\end{gather*}
 We note by $\gamma^{\sharp}_{i,k} = \gamma_{i,k}$ or $\gamma^{*}_{i,k}$ to avoid additional terms due to the symmetry. Clearly, $C^{A}_{r}$ is continuous from ${\rm CO}_{\infty-}\times {\rm CO}_{\infty-}$ into ${\rm CO}_{\infty-}$. We can f\/ind $k_{0}, C_{0} >0$ large enough and $K> 0$, by referring to the proof of Theorem~\ref{theorem4}, such that
\begin{gather*}
\big\|C^{A}_{r}\big(I^{S}_{m}(F),I^{S}_{m}(G)\big)\big\|_{k,C}\leq K\big\|I^{S}_{m}(F)\big\|_{k_{0},C_{0}}\big\|I^{S}_{m}(G)\big\|_{k_{0},C_{0}}<\infty.
\end{gather*}

\begin{definition}\label{definition10}
We put
\begin{gather}\label{eq10}
I^{S}_{m}(F)\star^{A}_{\hbar}I^{S}_{m}(G)= I^{S}_{m}(F).I^{S}_{m}(G)
+ \sum_{r\geq 1}\frac{\hbar^{r}}{r!}C^{A}_{r}\big(I^{S}_{m}(F),I^{S}_{m}(G)\big).
\end{gather}
\eqref{eq10} def\/ines  a  deformation quantization of $\{\cdot,\cdot \}$ in the sense of Def\/inition~\ref{definition6}.
\end{definition}

Finally,  we have
\begin{proposition}\label{proposition5}
$\star^{A}_{\hbar}$ and $\star_{\hbar}$ are equivalent on the Connes space ${\rm CO}_{\infty-}$.
\end{proposition}

\begin{proof}
We put as in~\cite{GDi}, for all $I^{S}_{m}(F)\in {\rm CO}_{\infty-}$
\begin{gather*}
T_{1}I^{S}_{m}(F)=-\sum_{i,k\geq 1} \alpha_{k}D_{\gamma_{i,k}}D_{\gamma^{*}_{i,k}}I^{S}_{m}(F).
\end{gather*}
Then, integrating by parts  $\forall\, k,C>0$, there exist $k_{1}>0$ and $C_{1}>0$ large enough such that
\begin{gather*}
\big\|T_{1}I^{S}_{m}(F)\big\|_{k,C}\leq \sum (Ck^2+1)^{-1}\big\|I^{S}_{m}(F)\big\|_{k_{1},C_{1}}<\infty	.
\end{gather*}
So $T_{1}$ is continuous on ${\rm CO}_{\infty-}$. We put $T := \exp(\hbar T_{1})$. $T$ is the formal series of operators
\begin{gather*}
T I^{S}_{m}(F)=\sum_{r\geq 0}\frac{(-\hbar)^{r}}{r!}\sum_{\substack{i_{1},\ldots,i_{r}\geq 1\\k_{1},\ldots,k_{r}\geq 1}}\prod_{l=1}^{r}\alpha_{k_{l}}D_{\gamma_{i_{1},k_{1}}}\cdots D_{\gamma_{i_{r},k_{r}}}D_{\gamma^{*}_{i_{1},k_{1}}}\cdots D_{\gamma^{*}_{i_{r},k_{r}}}I^{S}_{m}(F),
\end{gather*}
where we have seen that
\begin{gather*}
T^{r}(I^{S}_{m}(F)) := \sum_{\substack{i_{1},\ldots,i_{r}\geq 1\\k_{1},\ldots,k_{r}\geq 1}}\prod_{l=1}^{r}\alpha_{k_{l}}
D_{\gamma_{i_{1},k_{1}}}\cdots D_{\gamma_{i_{r},k_{r}}}D_{\gamma^{*}_{i_{1},k_{1}}}\cdots D_{\gamma^{*}_{i_{r},k_{r}}}I^{S}_{m}(F).
\end{gather*}
Since
$
T^{r} = \underbrace{T_{1}\circ\cdots\circ T_{1}}_{r\textrm{-times}}$,
we get that $T^{r}$ is continuous and as a result $T$ is continuous.

We note by $\langle \cdot,\cdot \rangle_{c} : \mathcal{H} \times \mathcal{H}^{*}\to \mathbb{R}$ the canonical pairing between $\mathcal{H}$ and $\mathcal{H}^{*}$. Then, according to~\cite{GDi}, we have the formula
\begin{gather*}
\Phi_{\gamma_{1},\gamma^{*}_{1}}\star^{A}_{\hbar} \Phi_{\gamma_{2},\gamma^{*}_{2}}=\exp\big[\hbar\big(\langle(A+I)\gamma_{1},\gamma^{*}_{2}\rangle_{c}
 + \langle(A-I)\gamma_{2},\gamma^{*}_{1}\rangle_{c}\big)\big]\Phi_{\gamma_{1}+\gamma_{2},\gamma^{*}_{1}+\gamma^{*}_{2}}.
\end{gather*}
Then as in~\cite{GDi}, we get
\begin{gather*}
T\big(\Phi_{\gamma_{1},\gamma^{*}_{1}}\star^{A}_{\hbar} \Phi_{\gamma_{2},\gamma^{*}_{2}}\big)= T(\Phi_{\gamma_{1},\gamma^{*}_{1}})\star_{\hbar}T(\Phi_{\gamma_{2},\gamma^{*}_{2}}).
\end{gather*}
This proves the proposition since the Wick exponentials are dense in the Connes spa\-ce \linebreak ${\rm CO}_{\infty-}$.
\end{proof}

\begin{remark}
In the Connes space ${\rm CO}_{\infty-}$ endowed with Stratonovich chaos, unlike in~\cite{GDi}, Moyal star-product and normal star-product($A = I$) are obviously equivalent. We can suppose that equivalences of \cite{R-S} with the Hida test functional space remain true because the Connes spa\-ce~${\rm CO}_{\infty-}$ is very small.
\end{remark}

\pdfbookmark[1]{References}{ref}
 \LastPageEnding

\end{document}